\documentclass[10pt,reqno]{amsart}
\usepackage{mathrsfs}
\usepackage{amssymb,amsmath,amscd,amsthm,amsfonts}
\newcommand{\f}{\frac}

\newcommand{\bea}{\begin{eqnarray}}
\newcommand{\eea}{\end{eqnarray}}
\newcommand{\bna}{\begin{eqnarray*}}
\newcommand{\ena}{\end{eqnarray*}}

\renewcommand{\le}{\left}
\newcommand{\ri}{\right}
\newcommand{\s}{\sigma}

\newtheorem{Theorem}{Theorem}[section]
\newtheorem{Lemma}[Theorem]{Lemma}

\newtheorem{Conjecture}{Conjecture}[section]

\begin{document}

\bigskip
\centerline{\sc \large A Prime number theorem for Rankin-Selberg}
\centerline{\sc \large $L$-functions over number fields}
\bigskip
\centerline{ Tim Gillespie and Guanghua Ji}

\openup 0.9\jot

\bigskip

\bigskip

\centerline{\scshape Abstract}
\begin{quote}
{\footnotesize In this paper we define a Rankin-Selberg $L$-function
attached to automorphic cuspidal representations of $GL_m({\Bbb
A}_E)\times GL_{m'}({\Bbb A}_{F})$ over cyclic algebraic number
fields $E$ and $F$ which are invariant under the Galois action, by
exploiting a result proved by Arthur and Clozel, and prove a prime
number theorem for this $L$-function. We end by rewriting the asymptotic formulas as sums over primes with a restriction on the modular degrees of the field extensions $E$ and $F$.  
\\{2000 Mathematics Subject Classification: 11F70, 11M26, 11M41.}}
\end{quote}

\bigskip

\section{Introduction}

Let $E$ be a Galois extension of $\mathbb{Q}$ of degree $\ell$. Let
$\pi$ be an automorphic cuspidal representation of $GL_m({\Bbb
A}_E)$ with unitary central character. Then the finite part
$L$-function attached to $\pi$ is given by the product of local
factors for $Re s=\sigma>1$, $L(s, \pi )=\prod_v L_\nu(s, \pi)$ (see \cite{GJ})
 where $L_\nu(s,\pi)=\prod_{j=1}^m \Big(1-
\f{\alpha_{\pi}(v,j)}{q_{v}^{s}}\Big)^{-1}$ and
$\alpha_{\pi}(v,j), 1\leq j\leq m$ are complex numbers given by the
Langlands correspondence, and $q_\nu$ denotes the cardinality of the residue at the place $\nu$. If $\pi_v$ is
ramified, we can also write the local factors at ramified places $v$
in the same form (1.1) with the convention that some of the
$\alpha_{\pi}(v,j)$ may be zero.  
For two automorphic cuspidal representations $\pi$ and
$\pi'$ of $GL_m({\Bbb A}_E)$ and $GL_{m'}({\Bbb A}_E)$,
respectively, denote the usual Rankin-Selberg $L$-function by \begin{equation}
L(s, \pi \times \widetilde{ \pi}') =\prod_v L_\nu (s, \pi\times
\widetilde{ \pi}')= \prod_{j=1}^m\prod_{i=1}^{m'}
\le(1-\f{\alpha_\pi (v,j) \overline{\alpha_{\pi'}
(v,i)}}{q_v^{s}}\ri)^{-1}.
  \end{equation} 

For $\sigma>1$, we have
\bna
\frac{L'}{L}\left(s,\pi\times\widetilde{\pi}'\right) =
-\sum_{n=1}^{\infty}\frac{\Lambda(n)a_{\pi\times\widetilde{\pi}'}(n)}{n^s},
\ena
see $\S2$, for the detailed definition of
$a_{\pi\times\widetilde{\pi}'}(n)$.
By a prime number theorem for
Rankin-Selberg $L$-functions $L(s, \pi \times\widetilde{\pi}')$, we
mean the asymptotic behavior of the sum
\bea
\sum_{n\leq x}\Lambda(n)a_{\pi\times\widetilde{\pi}'}(n).
\eea

A prime number theorem for Rankin-Selberg $L$-functions with $\pi$
and $\pi'$ being classical holomorphic cusp forms has been studied
by several authors. Recently, Liu and Ye \cite{LiuYe4} computed a
revised version of Perron's formula. Using the new Perron's formula,
the authors proved a prime number theorem for
Rankin-Selberg $L$-functions over $\mathbb{Q}$ without assuming the
Generalized Ramanujan Conjecture. Following the method in
\cite{LiuYe4}, we obtain a prime number theorem of the
Rankin-Selberg $L$-functions defined over a number field $E$.

\begin{Theorem}
Let $E$ be Galois extension of $\mathbb{Q}$ of degree $\ell$. Let
$\pi$ and $\pi'$ be irreducible unitary cuspidal representations of
$GL_{m}({\Bbb A}_{\Bbb E})$ and $GL_{m'}({\Bbb A}_{\Bbb E})$,
respectively. Assume that at least one of $\pi$ or $\pi'$ is
self-contragredient. Then \bna &&\sum_{n\leq
x}\Lambda(n)a_{\pi\times\widetilde{\pi}'}(n) \nonumber
\\
&& = \le\{\begin{array}{l} \frac{\displaystyle
x^{1+i\tau_0}}{\displaystyle 1+i\tau_0} +O\{x\exp(-c\sqrt{\log x})\}
\\
\hspace{20mm} \text{if}\  \pi'\cong\pi\otimes|\det|^{i\tau_0}\
\text{for}\ \text{some}\ \tau_0\in{\Bbb R};
\\
O\{x\exp(-c\sqrt{\log x})\}
\\
\hspace{20mm} \text{if}\ \pi'\not\cong\pi\otimes|\det|^{it}\
\text{for}\  \text{any}\   t\in{\Bbb R}.
\end{array}
\right.
\ena
\end{Theorem}

\smallskip

Let $E$ be a cyclic Galois extension of $\mathbb{Q}$ of degree
$\ell$. Let $\pi$ be an automorphic cuspidal representation of
$GL_m({\Bbb A}_E)$ with unitary central character. Suppose that
$\pi$ is stable under the action of $\textrm{Gal}(E/\mathbb{Q})$.
Thanks to Arthur and Clozel \cite{AC}, $\pi$ is the base change lift of
exactly $\ell$ nonequivalent cuspidal representations
\begin{equation}  \pi_{\mathbb{Q}},
\pi_{\mathbb{Q}}\otimes\eta_{E/\mathbb{Q}}, ...,
\pi_{\mathbb{Q}}\otimes\eta_{E/\mathbb{Q}}^{\ell-1} \end{equation} of $GL_m(\Bbb
A_{\Bbb Q})$, where $\eta_{E/\mathbb{Q}}$ is a nontrivial
character of $\Bbb A_{\Bbb Q}^{\times}/{\Bbb Q}^{\times}$ attached
to the field extension $E$ according to class field theory.
Consequently, we have $L(s, \pi)=L(s, \pi_{\mathbb{Q}}) L(s,
\pi_{\mathbb{Q}}\otimes\eta_{E/\mathbb{Q}})\cdots
L(s,\pi_{\mathbb{Q}}\otimes\eta_{E/\mathbb{Q}}^{\ell-1})$
where
the $L$-functions on the right side are distinct.

Similarly, let $F$ be a cyclic Galois extension of $\mathbb{Q}$ of
degree $q$. Let $\pi'$ be an automorphic cuspidal representation of
$GL_{m'}({\Bbb A}_F)$ with unitary central character, and
suppose that $\pi'$ is stable under the action of
$\textrm{Gal}(F/\mathbb{Q})$.  Then we can write 
\begin{equation}L(s,\pi')=\prod_{j=0}^{q-1}L(s,\pi'_{\mathbb{Q}}\otimes\psi_{F/\mathbb{Q}}^j)    \end{equation} 
where $\pi'_{\mathbb{Q}}$ is an irreducible cuspidal
representation of $GL_{m'}({\Bbb A}_\mathbb{Q})$ and
$\psi_{F/\mathbb{Q}}$ is a nontrivial character of $\Bbb A_{\Bbb
Q}^{\times}/{\Bbb Q}^{\times}$ attached to the field extension $F$.
Then we define the Rankin-Selberg $L$-function over the different
number fields $E$ and $F$ by \bea
L(s,\pi\times_{BC}\widetilde{\pi}') = \prod_{{0\leq i\leq
\ell-1}\atop{0\leq j\leq
q-1}}L(s,\pi_{\mathbb{Q}}\otimes\eta_{E/\mathbb{Q}}^{i}\times
\widetilde{\pi'_{\mathbb{Q}}\otimes\psi_{F/\mathbb{Q}}^{j}}), \eea
where
$L(s,\pi\otimes\eta^{i}\times\widetilde{\pi'\otimes\psi^{j}})$,
$0\leq i\leq \ell-1,\;0\leq j\leq q-1$ are the usual Rankin-Selberg
$L$-functions over $\mathbb{Q}$ with unitary central characters.
Then for $\s>1$, we have \bna 
-\f{d}{ds}\log L(s,\pi\times_{BC}\widetilde{\pi}')= 
\sum_{n=1}^{\infty}\f{\Lambda(n)a_{\pi\times_{BC}\widetilde{\pi}'}(n)}{n^s},
\ena where
$$a_{\pi\times_{BC}\widetilde{\pi}'}(n) =\sum_{0\leq i\leq
\ell-1}\sum_{0\leq j\leq
q-1}a_{\pi_{\mathbb{Q}}\otimes\eta_{E/\mathbb{Q}}^{i}}
(n)a_{\widetilde{\pi'_{\mathbb{Q}}\otimes\psi_{F/\mathbb{Q}}^{j}}}(n).$$
By a prime number theorem for Rankin-Selberg $L$-functions $L(s, \pi
\times_{BC}\pi')$ over number fields $E$ and $F$, we mean the
asymptotic behavior of the sum \bea \sum_{n\leq
x}\Lambda(n)a_{\pi\times_{BC}\widetilde{\pi}'}(n) = \sum_{n\leq x}\sum_{0\leq
i\leq \ell-1}\sum_{0\leq j\leq q-1}\Lambda(n)a_{\pi_{\mathbb{Q}}
\otimes\eta_{E/\mathbb{Q}}^{i}}(n)
a_{\widetilde{\pi'_{\mathbb{Q}}\otimes\psi_{F/\mathbb{Q}}^{j}}}(n).
\eea

Using the main theorem in Liu and Ye \cite{LiuYe4}, we obtain a
prime number theorem over different number fields $E$ and $F$.

\begin{Theorem}
Let $E$ and $F$ be two cyclic Galois extensions of $\mathbb{Q}$ of
prime degrees $\ell$ and $q$, respectively, with $(\ell,q)=1$. Let $\pi$
and $\pi'$ be unitary automorphic cuspidal representations of
$GL_{m}({\Bbb A}_{E})$ and $GL_{m'}({\Bbb A}_{ F})$,
respectively. Assume that we have base change lifts as in  (1.4) and (1.5),
and suppose that $\pi_{\mathbb{Q}}$ is self contragredient, then
\bna
&&\sum_{n\leq x}\Lambda(n)a_{\pi\times_{BC}\widetilde{\pi}'}(n) \nonumber
\\
&& = \le\{\begin{array}{l}\frac{\displaystyle
x^{1+i\tau_0}}{\displaystyle 1+i\tau_0} +O\{x\exp(-c\sqrt{\log x})\}
\\
\hspace{4mm} \text{if}\
\pi'_{\mathbb{Q}}\otimes\psi_{F/\mathbb{Q}}^{j_0}
\cong\pi_{\mathbb{Q}}\otimes\eta_{E/\mathbb{Q}}^{i_0}\otimes|\det|^{i\tau_0}
\ \text{for}\ \text{some}\  \tau_0\in{\Bbb R}\text{  and some  }
i_0, j_0; 
\\
O\{x\exp(-c\sqrt{\log x})\}
\\
\hspace{4mm} \text{if}\
\pi'_{\mathbb{Q}}\otimes\psi_{F/\mathbb{Q}}^{j}
\ncong\pi_{\mathbb{Q}}\otimes\eta_{E/\mathbb{Q}}^{i}\otimes|\det|^{i\tau}\
\text{for}\  \text{any}\  i,j \text{   and   } \tau\in{\Bbb R}.
\end{array}
\right. \ena
\end{Theorem}

We end by rewriting Theorems 1.1 and 1.2 as sums over primes using Conjecture 2.1 and Hypothesis H to show that the main term comes from those primes which split completely in the field extension
\begin{Theorem}(1) Let the notations be as in Theorem 1.1.  Assume Hypothesis H and Conjecture 2.1 , then 
\begin{equation}\sum_{{p\leq x}\atop {\text{$p$ splits completely}}}(\log p)a_{\pi\times\widetilde{\pi}'}(p)=\frac{x^{1+i\tau_0}}{1+i\tau_0}+O\{x\exp(-c\sqrt{\log x})\}\text{ for }\pi\cong\pi'\otimes|\det|^{i\tau_0}  \nonumber  \end{equation} 

(2)  Let the notations be as in Theorem 1.2 and suppose for some $i_0,j_0$ and $\tau_0\in \mathbb{R}$ that $\pi_{\mathbb{Q}}\otimes\eta_{E/\mathbb{Q}}^{i_0}\cong\pi_{\mathbb{Q}}'\otimes\psi_{F/\mathbb{Q}}^{j_0}\otimes|\det|^{i\tau_0}$.  Assume Hypothesis H and Conjecture 2.1, also suppose that for any prime $p$ for which both $\pi_\nu$ and $\pi'_\omega$ are unramified for any $\nu|p$, $\omega|p$ that we have the following:  suppose that there exist primes $\mathfrak{p}$ and $\mathfrak{q}$ in the ring of integers of $E$ and $F$ ,respectively, lying above $p$ which also lie below  primes $\mathfrak{P}$ and $\mathfrak{Q}$ inside the ring of integers of $EF$ with the restriction $f_{\mathfrak{Q}/\mathfrak{q}}\leq f_p$ and $f_{\mathfrak{P}/\mathfrak{p}}\leq f_p$.  Then 
\begin{equation}\sum_{{p\leq x}\atop{\text{$p$ splits completely in $EF$}}}(\log p)a_{\pi\times_{BC}\widetilde{\pi}'}(p)=\frac{x^{1+i\tau_0}}{1+i\tau_0}+O\{x \exp(-c\sqrt{\log x})\} \nonumber  \end{equation} \nonumber  \end{Theorem}    
\medskip
{\it Remark}:   Note Theorem 1.3 says to obtain the main term we need only consider those summands for which $f_p=f_p'=1$.  Such conditions are useful in controlling sums over primes in the computation of the n-level correlation function attached to a cuspidal representation of $GL_n(\mathbb{A}_E)$ over a number field $E$ (see \cite{LiuYe1}). 
 \medskip

\smallskip

\section{Rankin-Selberg L-functions}
\setcounter{equation}{0}

In the section we recall some fundamental analytic properties of
Rankin-Selberg $L$-functions: absolute convergence of the Euler product, location of poles, and zero-free region.

Let $E$ be a Galois extension of $\mathbb{Q}$ of degree $\ell$. For
any prime $p$, we have
$E\otimes_{\mathbb{Q}}\mathbb{Q}_p=\oplus_{v|p}E_v$, where $v$
denotes a place of $E$ lying above $p$. Since $E$ is Galois over
$\mathbb{Q}$, all $E_v$ with $v|p$ are isomorphic. Denote by
$\ell_p$ the degree $[E_{\nu}:\mathbb{Q}_p]$ , by $e_p$ the order of
ramification, and by $f_p$ the modular degree of $E_v$ over
$\mathbb{Q}_p$ for $v|p$. Then we have $\ell_p=e_pf_p$, and
$q_v=p^{f_p}$ is the cardinality of the residue class field.
Let $\pi$ be an irreducible cuspidal representation of $GL_m({\Bbb
A}_E)$ with unitary central character.  Let $\pi'$ be an automorphic cuspidal representation of
$GL_{m'}(\mathbb{A}_E)$ with unitary central character. The finite-part
Rankin-Selberg $L$-function $L(s, \pi\times\widetilde{\pi}')$ is
given by the product of local factors. 
\begin{equation}L(s,\pi\times\tilde{\pi}')=\prod_{\nu<\infty}L_{\nu}(s,\pi\times\tilde{\pi}')  \nonumber \end{equation} and we denote
\bna L_p(s,\pi\times\widetilde{\pi}')
&=&\prod_{v|p}L_v(s,\pi_v \times \widetilde{ \pi}'_v)
=\prod_{v|p}\prod_{j=1}^m\prod_{i=1}^{m'} \Big(1-\f{\alpha_\pi
(v,j) \overline{\alpha_{\pi'} (v,i)}}{p^{f_{p}s}}\Big)^{-1}. \ena Then for
$\sigma>1$, we have \bea
\frac{L'}{L}(s,\pi\times\widetilde{\pi}') &=&
-\sum_{v}\sum_{j=1}^{m}\sum_{i=1}^{m'}\sum_{k\geq 1}\frac{f_p\log
p}{p^{kf_{p}s}}
\alpha^k_\pi(v,j) \overline{\alpha^k_{\pi'} (v,i)}\nonumber\\
&=&
-\sum_{n=1}^{\infty}\frac{\Lambda(n)a_{\pi\times\widetilde{\pi}'}(n)}{n^s},
\eea
where
$$a_{\pi\times\widetilde{\pi}'}(p^{kf_p})=
\sum_{\nu|p}f_p\Big(\sum_{j=1}^{m}\alpha_\pi(\nu,j)^k\Big)
\Big(\sum_{i=1}^{m'}\overline{\alpha_{\pi'}(\nu,i)}^k\Big),$$ and
$a_{\pi\times\widetilde{\pi}'}(p^{k})=0$, if $f_p\nmid k$.

We will use the Rankin-Selberg $L$-functions $L(s, \pi \times
\widetilde\pi')$ as developed by Jacquet, Piatetski-Shapiro, and
Shalika \cite{JacPiaSha}, Shahidi \cite{Sha1}, and Moeglin and
Waldspurger \cite{MoeWal}.  We will need the following properties of $L(s,\pi\times\widetilde{\pi}')$ 

\medskip

{\bf RS1}. The Euler product for $L(s,\pi\times\widetilde{\pi}')$ in
(2.1) converges absolutely for $\s>1$ (Jacquet and Shalika
\cite{JacSha1}).

\medskip

{\bf RS2}. Denote $\alpha(g)=|\det(g)|$. When $\pi'\not\cong
\pi\otimes\alpha^{it}$ for any $t\in{\Bbb R}$,
$L(s,\pi\times\widetilde{\pi}')$ is holomorphic. When $m=m'$ and
$\pi' \cong \pi\otimes \alpha^{i\tau_0} $ for some $\tau_0\in\Bbb
R$, the only poles of $L(s, \pi \times \widetilde\pi ')$ are simple
poles at $s=i\tau_0$ and $1+i\tau_0$ (Jacquet and Shalika \cite{JacSha1},
Moeglin and Waldspurger \cite{MoeWal}).

\medskip
Finally, we note the reason for the self-contragredient assumption is that one must apply the following zero-free region due to Moreno to obtain the error term as in the Theorems
\medskip

{\bf RS3}.
$L(s,\pi\times\widetilde{\pi}')$ is non-zero in $\s\ge 1$ (Shahidi
\cite{Sha1}).  Furthermore, if at least one of $\pi$ or $\pi'$ is
self-contragredient, it is zero-free in the region
\begin{equation}
\sigma > 1-\frac{c}{\log(Q_{\pi}Q_{\pi'}(|t|+2))}, \quad |t|\geq 1
\end{equation}
where $c$ is an explicit constant depending only on $m$ and $n$ (see
Sarnak \cite{Sa}, Moreno \cite{Mor} or Gelbert, Lapid and Sarnak
\cite{GLS}).

\medskip

Let $L(s,\pi\times_{BC}\widetilde{\pi}')$ be the Rankin-Selberg
$L$-function over the number fields $E$ and $F$,
$$
L(s,\pi\times_{BC}\widetilde{\pi}')=\prod_{{0\leq i\leq
\ell-1}\atop{0\leq j\leq
q-1}}L(s,\pi_{\mathbb{Q}}\otimes\eta_{E/\mathbb{Q}}^{i}
\times\widetilde{\pi'_{\mathbb{Q}}\otimes\psi^j_{F/\mathbb{Q}}})
$$
where $L(s,\pi_{\mathbb{Q}}\otimes\eta_{E/\mathbb{Q}}^{i}
\times\widetilde{\pi'_{\mathbb{Q}}\otimes\psi^j_{F/\mathbb{Q}}})$
is the usual Rankin-Selberg $L$-function on $GL_{m}\times GL_{n}$
over $\Bbb Q$. Hence
$L(s,\pi\times_{BC}\widetilde{\pi}')$ will have similar analytic
properties as the usual Rankin-Selberg $L$-functions.  We will need the following bound for the local parameters proved in \cite{RudSa}  
\begin{equation} |\alpha_{\pi}(j,\nu)|\leq p^{f_p(1/2-1/(m^2\ell+1))}\text{  for  }\nu|p   \end{equation}
this holds for both $\pi_\nu$ ramified and unramified.  When $\pi_\nu$ is unramified the generalized Ramanujan conjecture claims that $|\alpha_{\pi}(j,\nu)|=1$.
The best known bound toward this conjecture over an arbitrary number field is $\alpha_\pi(j,\nu)\leq p^{{f_p}/9}$ for $m=2$ \cite{KimSha}.  We will not assume the generalized Ramanujan conjecture, but assume a bound $\theta_p$ toward it for any $p$ which is unramified and does not split completely in $E$.    
\begin{Conjecture}For any $p$ which is unramified and does not split completely in E, we have for any $\nu|p$ that 
\begin{equation}  |\alpha_\pi(j,\nu)|\leq p^{f_p\theta_p}  \nonumber  \end{equation}
where $\theta_p=1/2-1/(2f_p)-\epsilon$ for a small $\epsilon>0$.  
  \nonumber  \end{Conjecture}  
Note that for $\pi_\nu$ unramified we have $e_p=1$ and hence $f_p=\ell_p$ where $\ell_p=[E_\nu:\mathbb{Q}_p]$.  Since $p$ does not split completely in $E$, we know that $f_p\geq 2$.  Thus Conjecture 2.1 is known for $m=2$ according to (2.2).  It is trivial for $m=1$.  Since $f_p|\ell$, Conjecture 2.1 is known when all prime factors of $\ell$ are $>(m^2+1)/2$.  For $m=3$ this means that any $p|\ell$ is $\geq 7$, while for $m=4$, Conjecture 2.1 is true when any $p|\ell$ is $\geq 11$.  We end this section by recalling Hypothesis H from \cite{LiuYe1}
\smallskip   

$\mathbf{Hypothesis}$ $\mathbf{H}$  {\it Let $\pi$ be an automorphic cuspidal representation of $GL_m(\mathbb{A}_{E})$ with unitary central character, then for any fixed $k\geq 2$}  
\begin{equation}  \sum_{p}\frac{\log^2 p}{p^{kf_p}}\sum_{\nu|p}\Big|\sum_{1\leq j\leq m}\alpha_{\pi}(i,\nu)^k\Big|^2<\infty  \nonumber \end{equation}

\smallskip

\section{Proof of Theorem 1.1 }
\setcounter{equation}{0}

Let $\pi$ and $\pi^{\prime}$ be as in the Theorem 1.1, we will first
need a modified version of Lemma 4.1 of \cite{LiuYe4}. It is a
weighted prime number theorem in the diagonal case.  With the same modifications made as in \cite{LiuYe1} Lemma 6.1 over a number field the proof follows as in \cite{LiuYe4}.   

\begin{Lemma}
Let $\pi$ be a self-contragredient automorphic irreducible cuspidal
representation of $GL_m$ over E.  Then \bna \sum_{n\leq
x}\left(1-\frac{n}{x}\right)\Lambda(n)a_{\pi\times\widetilde{\pi}}(n)
=\frac{x}{2}+O\lbrace x\exp(-c\sqrt{\log x})\rbrace \ena \end{Lemma}

The next lemma again closely follows \cite{LiuYe4}, and allows the
removal of the weight $(1-\frac{n}{x})$ from the previous lemma.
The proof involves a standard argument due to de la Vallee Poussin.
\begin{Lemma}
Let $\pi$ be a self-contragredient automorphic irreducible cuspidal
representation of $GL_m$ over E.  Then
\begin{equation}
\sum_{n\leq x}\Lambda(n)a_{\pi\times\widetilde{\pi}}(n)=x+O\lbrace
x\exp(-c\sqrt{\log x})\rbrace.\end{equation}
\end{Lemma}
$\emph{Proof.}$ Since the coefficients of the left hand side of
(3.7) are non-negative, the proof follows as in Lemma 5.1 of
\cite{LiuYe4} with no modification. $\quad\square$

The next lemma also follows exactly as in Lemma 5.2 of
\cite{LiuYe4}, and doesn't require $\pi$ to be self-contragredient.
The proof is an application of a Tauberian theorem of Ikehara.
\begin{Lemma}For any automorphic irreducible cuspidal unitary representation
$\pi$ of $GL_m$ over the number filed E, we have
\begin{equation}  \sum_{n\leq x} \Lambda(n)a_{\pi\times\widetilde{\pi}}(n)\thicksim x.
\end{equation}
\end{Lemma}

\emph{Proof of Theorem 1.1.}  We suppose throughout that $\pi$ is
self-contragredient. When $\pi'\cong\pi$, the theorem reduces to
Lemma 3.1.  We will first consider the case when $\pi$ and $\pi'$
are twisted equivalent, so suppose that
$\pi'\cong\pi\otimes|\textrm{det}|^{i\tau_0}$ for some
$\tau_{0}\in\mathbb{R}$.  By Lemma 3.2, we obtain a bound for the
short sum
\begin{equation}
\sum_{x<n\leq x+y}\Lambda(n)a_{\pi\times\widetilde{\pi}'}(n)\ll y
\nonumber
\end{equation}
for $y\gg x\exp(-c\sqrt{\log x})$.  $\pi'$ is not necessarily
self-contragredient; nevertheless, by Lemma 3.3, we get for $0<y\leq
x$ that
\begin{equation}
\sum_{x<n\leq x+y } \Lambda(n)a_{\pi\times\widetilde{\pi}'}(n)\ll
\sum_{x<n\leq 2x}\Lambda(n)a_{\pi\times\widetilde{\pi}'}(n)\ll x.
\nonumber
\end{equation}
By definition of the coefficients $a_{\pi\times\widetilde{\pi}'}(n)$, we
have that for $n=p^{kf_p}$
\bna
|\Lambda(n)a_{\pi\times\widetilde{\pi}'}(n)|
&\leq&
\log p\sum_{\nu|p}f_p\Big|\sum_{j=1}^{m}\alpha_\pi(\nu,j)^k\Big|
\Big|\sum_{i=1}^{m'}\bar{\alpha}_{\pi'}(\nu,i)^k\Big|  \nonumber\\
&\leq& \log
p\Big(\sum_{\nu|p}f_p\Big|\sum_{j=1}^{m}\alpha_\pi(\nu,j)^k\Big|^2\Big)^{1/2}
\Big(\sum_{\nu|p}f_p\Big|\sum_{i=1}^{m'}\bar{\alpha}_{\pi'}(\nu,i)^k\Big|^{2}\Big)^{1/2}
\nonumber\\
&=&
(\log p)a_{\pi\times\widetilde{\pi}}(n)^{1/2}a_{\pi'\times\widetilde{\pi}'}(n)^{1/2}.
\nonumber\ena
Thus we have,
\bea \sum_{x<n\leq x+y}|\Lambda(n)a_{\pi\times\widetilde{\pi}'}(n)|
&\leq&
\sum_{x<n\leq x+y}\Lambda(n)a_{\pi\times\widetilde{\pi}}(n)^{1/2}a_{\pi'\times\widetilde{\pi}'}(n)^{1/2}
\nonumber\\
&\leq&
\Big(\sum_{x<n\leq x+y}\Lambda(n)a_{\pi\times\widetilde{\pi}}(n)\Big)^{1/2}
\Big(\sum_{x<n\leq x+y}\Lambda(n)a_{\pi'\times\widetilde{\pi}'}(n)\Big)^{1/2}\nonumber\\
&\ll& \sqrt{yx}.
\eea
now the rest of the proof follows exactly as in \cite{LiuYe4} $\square$
\smallskip

\section{Proof of Theorem 1.2}
\setcounter{equation}{0}

Let $E$ and $F$ be two cyclic Galois extensions of $\Bbb Q$ of
degree $\ell$ and $q$, respectively. Let $\pi$ and $\pi'$ be
irreducible unitary cuspidal representations of $GL_m({\Bbb A}_E)$
and $GL_{m'}({\Bbb A}_F)$ with unitary central characters. In
section 2, we denote by
$$
L(s,\pi\times_{BC}\widetilde{\pi}')
=\prod_{{0\leq i\leq\ell-1}\atop{0\leq j\leq q-1}}
L(s,\pi_{\mathbb{Q}}\otimes\eta_{E/\mathbb{Q}}^{i}\times
\widetilde{\pi'_{\mathbb{Q}}\otimes\psi_{F/\mathbb{Q}}^{j}}),
$$
where
$L(s,\pi\otimes\eta_{E/\mathbb{Q}}^{i}\times\widetilde{\pi'\otimes\psi_{F/\mathbb{Q}}^{j}})$,
$0\leq i\leq\ell-1,\;0\leq j\leq q-1$ are the usual Rankin-Selberg
$L$-functions over $\mathbb{Q}$ with unitary central characters.

\begin{Lemma}
Suppose that
$\pi_{\mathbb{Q}}\otimes\eta_{E/\mathbb{Q}}^{i_0}\cong\pi'_{\mathbb{Q}}
\otimes\psi_{F/\mathbb{Q}}^{j_0}\otimes|\det|^{i\tau_0},$ for some
$0\leq i_0\leq\ell-1$, $0\leq j_0\leq q-1$ and $\tau_0\in
\mathbb{R}$. Then
\begin{equation} \pi_{\mathbb{Q}}\otimes\eta_{E/\mathbb{Q}}^{i_0}
\cong\pi'_{\mathbb{Q}}\otimes\psi_{F/\mathbb{Q}}^{j}\otimes |\det|^{i\tau} \nonumber
\end{equation}
implies that $\tau=\tau_{0}$ and $j=j_{0}$. Moreover, if
\begin{equation}
\pi_{\mathbb{Q}}\otimes\eta_{E/\mathbb{Q}}^{i}\cong\pi'_{\mathbb{Q}}\otimes\psi_{F/\mathbb{Q}}^{j}\otimes
|\det|^{i\tau} \nonumber \end{equation} for some $i$ and $j$, and
$\tau\in\mathbb{R}$, then $\tau=\tau_{0}$.
\end{Lemma}

\textit{Proof}. By class field theory, $\eta_{E/\mathbb{Q}},
\psi_{F/\mathbb{Q}}$ are finite order idele class characters, so
they are actually primitive Dirichlet characters. Assume that
$$
\pi_{\mathbb{Q}}\otimes\eta_{E/\mathbb{Q}}^{i_0}
\cong\pi'_{\mathbb{Q}}\otimes\psi_{F/\mathbb{Q}}^{j}\otimes|\det|^{i\tau},
$$
for some $0\leq i\leq \ell-1$, $0\leq j\leq q-1$ and
$\tau\in\mathbb{R}$. Then we have
$$
\pi'_{\mathbb{Q}}\otimes\psi_{F/\mathbb{Q}}^{j_0}\otimes|\det|^{i\tau_0}
\cong
\pi'_{\mathbb{Q}}\otimes\psi_{F/\mathbb{Q}}^{j}\otimes|\det|^{i\tau}.$$
For any unramified $p$, we get
$$\{\alpha_{\pi'_{\mathbb{Q}}}(p,j)
\psi_{F/\mathbb{Q}}^{j_0}{|p|}_p^{i\tau_0}\}_{j=1}^{m}
=\{\alpha_{\pi'_{\mathbb{Q}}}(p,j)\psi_{F/\mathbb{Q}}^{j}(p){|p|}_p^{i\tau}\}_{j=1}^{m}.$$
Hence,
$$(\psi_{F/\mathbb{Q}}^{j_0}p^{-i\tau_0})^m
=(\psi_{F/\mathbb{Q}}^{j}(p)p^{-i\tau})^m.$$ Since
$\psi_{F/\mathbb{Q}}$ is of finite order, we get by multiplicity one for characters $\tau=\tau_0$, so
that $j=j_{0}$.  The last conclusion of the lemma follows from the
same argument just given. $\quad\square$

\begin{Lemma} Suppose that $\pi_{\mathbb{Q}}\otimes
\eta_{E/\mathbb{Q}}^{i_0}\cong
\pi_{\mathbb{Q}}\otimes\psi_{F/\mathbb{Q}}^{j_0}\otimes|\det|^{i\tau_0}$
for some $0\leq i_0 \leq \ell-1$, $0\leq j_0 \leq q-1$ and
$\tau_0\in \mathbb{R}$. Then the number of twisted equivalent pairs
$( \pi_{\mathbb{Q}}\otimes\eta_{E/\mathbb{Q}}^i,
\pi'_{\mathbb{Q}}\otimes\psi_{F/\mathbb{Q}}^j)$ with $0\leq i\leq
\ell-1$, $0\leq j\leq q-1$ divides the greatest common divisor of
$\ell$ and $q$.
\end{Lemma}

{\it Proof.} By relabeling the collection
$\{\pi_{\mathbb{Q}}\otimes\eta_{E/\mathbb{Q}}^i\}_{0\leq i \leq
\ell-1}$ if necessary we may assume that
$\pi_{\mathbb{Q}}\cong\pi'_{\mathbb{Q}}\otimes\psi_{F/\mathbb{Q}}^{j_0}\otimes|\det|^{i\tau_0}$.
Now let $G=(\{\pi_{\mathbb{Q}}\otimes\eta_{E/\mathbb{Q}}^i,0\leq
i\leq\ell-1\},*)$ where we define
\begin{equation}
\pi_{\mathbb{Q}}\otimes\eta_{E/\mathbb{Q}}^{i_1}*\pi_{\mathbb{Q}}
\otimes\eta_{E/\mathbb{Q}}^{i_2}=\pi_{\mathbb{Q}}\otimes\eta_{E/\mathbb{Q}}^{i_1+i_2}.
\nonumber
\end{equation}
Since the character $\eta_{E/\mathbb{Q}}$
has order $\ell$ we have $G\cong \mathbb{Z}/\ell\mathbb{Z}$.  Now
let
\begin{equation}
H=\{\pi_{\mathbb{Q}}\otimes\eta_{E/\mathbb{Q}}^i: \exists\, 0\leq j
\leq q-1, \tau\in\mathbb{R}, \textrm{ such that }
\pi_{\mathbb{Q}}\otimes\eta_{E/\mathbb{Q}}^i\cong\pi'_{\mathbb{Q}}\otimes
\psi_{F/\mathbb{Q}}^j\otimes|\det|^{i\tau} \}. \nonumber
\end{equation} By hypothesis, we have $\pi_{\mathbb{Q}}\in H$. Assume that
$$
\pi_{\mathbb{Q}}\otimes\eta_{E/\mathbb{Q}}^{i_1}
\cong\pi'_{\mathbb{Q}}\otimes\psi_{F/\mathbb{Q}}^{j_1}\otimes|\det|^{i\tau_1}
$$ and
$$\pi_{\mathbb{Q}}\otimes\eta_{E/\mathbb{Q}}^{i_2}\cong\pi'_{\mathbb{Q}}
\otimes\psi_{F/\mathbb{Q}}^{j_2}\otimes|\det|^{i\tau_2},
$$
then
\begin{eqnarray}
\pi_{\mathbb{Q}}\otimes\eta_{E/\mathbb{Q}}^{i_1-i_2}
&\cong&
\pi'_{\mathbb{Q}}\otimes \psi_{F/\mathbb{Q}}^{j_1}\otimes\eta_{E/\mathbb{Q}}^{-i_2}\otimes|\det|^{i\tau_1}
\nonumber \\
&\cong&
\pi_{\mathbb{Q}}\otimes\psi_{F/\mathbb{Q}}^{j_1-j_2}\otimes|\det|^{i(\tau_1-\tau_2)}  \nonumber  \\
&\cong&
\pi'_{\mathbb{Q}}\otimes\psi_{F/\mathbb{Q}}^{j_0+j_1-j_2}\otimes|\det|^{i(\tau_1-\tau_2+\tau_0)}.
\nonumber
\end{eqnarray} Hence $H$ is a subgroup of $G$.
By Lemma 4.1, each $\pi_{\mathbb{Q}}\otimes \eta_{E/\mathbb{Q}}^i$
is twisted equivalent to at most one
$\pi'_{\mathbb{Q}}\otimes\psi_{F/\mathbb{Q}}^j$ so by Lagrange's
theorem the number of twisted equivalent pairs divides $\ell$, and
by symmetry of the above argument we also have that it divides $q$,
so the lemma follows. $\quad\square$

The above lemmas are simple but give the following: the second
conclusion says that we can have at most one twisted equivalent
pair when $\ell$ and q are relatively prime, and the first conclusion
says that if the $L$-function $L(s,\pi\times_{BC}\pi')$ has poles at
$1+i\tau_{0}$ and $i\tau_{0}$, then these are the only poles, with
orders possibly bigger than one. If one considers the diagonal case
\bea L(s,\pi\times_{BC}\widetilde{\pi}) &=&
\prod_{i=0}^{\ell-1}\prod_{j=0}^{q-1}
L(s,\pi_{\mathbb{Q}}\otimes\eta_{E/\mathbb{Q}}^{i}\times
\widetilde{\pi_{\mathbb{Q}}\otimes\eta_{E/\mathbb{Q}}^{j}})\nonumber\eea then we get by $\mathbf{RS2}$ a simple pole of order
$\ell$ at $s=1$ for each factor in the left hand side, since the
other factors on the right hand side are nonzero at $s=1$ by
$\mathbf{RS3}$, so that this differs by the classical case in that
we get multiple poles.  Now assuming $\pi_{\mathbb{Q}}$ to be
self-contragredient, and applying Theorem 1.1 to the $L$-function
$L(s,\pi\times_{BC}\widetilde{\pi}')$, we can use the zero-free
region in $\mathbf{RS3}$ to obtain the same error term as before since \bna
L(s,\pi_{\mathbb{Q}}\otimes\eta_{E/\mathbb{Q}}^{i-1}\times\widetilde{\pi'_{\mathbb{Q}}\otimes
\psi_{F/\mathbb{Q}}^{j-1}})
=L(s,\pi_{\mathbb{Q}}\times\widetilde{\pi'_{\mathbb{Q}}\otimes
\psi_{F/\mathbb{Q}}^{j-1}\otimes\eta_{E/\mathbb{Q}}^{-(i-1)}}). \ena
We can apply the zero-free region to all the factors in the
definition of
 $L(s,\pi\times_{BC}\widetilde{\pi}')$, to get the same error term as in Theorem 1.1.
Thus Theorem 1.2 follows directly from Theorem 1.1.

\section { Sums over primes}
Note that Theorem 1.1 says  
\begin{equation}\sum_{\underset{p^{kf_p}\leq x}{p,k}}\log(p)\sum_{\nu|p}f_p\sum_{i=1}^{m}\sum_{j=1}^{m'}\alpha_{\pi}(i,\nu)^k\overline{\alpha_{\pi'}(j,\nu)}^k   \nonumber  \end{equation}
\begin{equation}=\frac{x^{1+i\tau_0}}{1+i\tau_0}+O\{x\exp(-c\sqrt{\log x})\} \text{  for  }\pi\cong\pi'\otimes|\det|^{i\tau_0}    \end{equation}

We can apply the bound in (2.3) to the sum 
\begin{equation}\sum_{\underset{p^{f_pk}\leq x}{k> (m^2\ell+1)/2}}\log(p)\sum_{\nu|p}f_p\sum_{i=1}^{m}\sum_{j=1}^{m'}\alpha_{\pi}(i,\nu)^k\overline{\alpha_{\pi'}(j,\nu)}^k  \nonumber  \end{equation}
 
\begin{equation} \ll\sum_{\underset{p^{kf_p\leq x}}{k>(m^2\ell+1)/2}}(\log p) p^{2kf_p(1/2-1/(m^2\ell+1))}\ll x\Big( \sum_{n\leq x}\frac{\Lambda(n)}{n^{1+\epsilon}}\Big)   \nonumber \end{equation}  
for small $\epsilon>0$.  By partial summation we get 
\begin{equation}  \sum_{n\leq x}\frac{\Lambda(n)}{n^{1+\epsilon}}=\Big(x+O\{x\exp(-c\sqrt{\log x})\} \Big)\frac{1}{x^{1+\epsilon}}-\int_{1}^{x}\Big(t+O\{t\exp(-c\sqrt{\log t})\}\Big)\frac{1}{t^{2+\epsilon}}(-1-\epsilon)dt  \nonumber  \end{equation}
\begin{equation}=O\{x\exp(-c\sqrt{\log x})\}   \nonumber \end{equation} 
Note that Hypothesis H gives that for fixed $k\geq 2$
\begin{equation} \sum_{p}\frac{|a_{\pi\times\tilde{\pi}'}(p^{kf_p})|(\log(p^{kf_p}))^2}{p^{kf_p}}<\infty  \nonumber \end{equation} 
using this and partial summation for fixed $k\geq 2$ we can write for $y=\exp(\exp(c\sqrt{\log x}))$ and $x$ sufficiently large (note that $m=m'$) 
\begin{equation}\sum_{p^{kf_p}\leq x}\log p \sum_{\nu|p}f_p\sum_{i=1}^{m}\sum_{j=1}^{m'}\alpha_{\pi}(i,\nu)^k\overline{\alpha_{\pi'}(j,\nu)}^k   \nonumber  \end{equation}
\begin{equation}\ll x \sum_{p^{kf_p}\leq y}\frac{\log^2 p}{p^{kf_p}}\sum_{\nu|p}\Big|\sum_{i=1}^{m}\alpha_{\pi}(i,\nu)^k\Big|^2\frac{1}{\log p}    \nonumber  \end{equation} 
\begin{equation}=x\Big(O(1)\left. \frac{1}{\log t}\right|_{2}^{y}-\int_{2}^{y}O(1)\frac{1}{t(\log t)^2}dt\Big)\ll x\exp\{-c\sqrt{\log x}\}  \nonumber \end{equation}

Finally, using Conjecture 2.1 we get 
\begin{equation}  \sum_{{p^{f_p}\leq x}\atop{\text{p not split}}}\log(p)a_{\pi\times\tilde{\pi}'}(p^{f_p})\ll\sum_{{p^{f_p} \leq x}\atop{\text{p not split}}}\log(p)f_p\sum_{\nu|p}\Big|\sum_{i=1}^{m}\alpha_{\pi}(i,\nu)\Big|^2    \nonumber \end{equation}   
\begin{equation}\ll \sum_{{p^{f_p} \leq x}\atop{\text{p not split}}}(\log p)p^{2f_p(1/2-1/(2f_p)-\epsilon)}  \ll x\Big(\sum_{n\leq x}\frac{\Lambda(n)}{n^{1+\epsilon}}\Big)\ll x\exp\{-c\sqrt{\log x}\}   \nonumber \end{equation}
  We can do a similar calculation for Theorem 1.2 by first noting that 
\begin{equation} \sum_{\nu|p}f_p\sum_{i=1}^{m}\sum_{j=1}^{m'}\alpha_{\pi}(i,\nu)^k=\sum_{a=0}^{\ell-1}\sum_{i=1}^{m}\alpha_{\pi_\mathbb{Q}\otimes\eta_{E/\mathbb{Q}}^{a}}(i,p)^{f_pk}   \end{equation} 
for $n=p^{kf_p}$, and similarly 
\begin{equation}\sum_{\omega|p}f_p'\sum_{j=1}^{m'}\alpha_{\pi'}(j,\nu)^k=\sum_{b=0}^{q-1}\sum_{j=1}^{m'}\alpha_{\pi_{\mathbb{Q}}\otimes\psi_{F/\mathbb{Q}}^b}(j,\nu)^{f_p'k}     \end{equation}
Thus we get 
\begin{equation} \sum_{n\leq x}\Lambda(n)a_{\pi\times_{BC}\widetilde{\pi'}}(n)=\sum_{\underset{}{p^{k_1f_p}=p^{k_2f_p'}\leq x}}(\log p)\sum_{\nu|p}\sum_{\omega|p}f_pf_p'\sum_{i=1}^{m}\sum_{j=1}^{m'}\alpha_{\pi}(i,\nu)^{k_1}\overline{\alpha_{\pi'}(j,\omega)}^{k_2}  \nonumber  \end{equation} 
again by (2.1)  
\begin{equation}\sum_{\underset{k_1>\min\{ 2/(m^2\ell+1),2/(m'^2q+1)\}}{p^{k_1f_p}=p^{k_2f_p'}\leq x}}(\log p)\sum_{\nu|p}\sum_{\omega|p}f_pf_p'\sum_{i=1}^{m}\sum_{j=1}^{m'}\alpha_{\pi}(i,\nu)^{k_1}\overline{\alpha_{\pi'}(j,\omega)}^{k_2}   \nonumber  \end{equation}
\begin{equation}\ll\sum_{\underset{k_1>\min\{ 2/(m^2\ell+1),2/(m'^2q+1)\}}{p^{k_1f_p}=p^{k_2f_p'}\leq x}}(\log p)p^{k_1f_p(1/2-1/(m^2\ell+1))+k_2f_p'(1/2-1/(m'^2q+1))}  \nonumber  \end{equation}
\begin{equation}\ll x\Big( \sum_{n\leq x}\frac{\Lambda(n)}{n^{1+\epsilon}}\Big)\ll x\exp(-c\sqrt{\log x})   \nonumber  \end{equation}
 Now consider the sum 
\begin{equation} \sum_{\underset{\text{ p not split in EF}}{p^{k_1f_p}=p^{k_2f_p'}\leq x}} (\log p)\sum_{\nu|p}\sum_{\omega|p}f_pf_p'\sum_{i=1}^{m}\sum_{j=1}^{m'}\alpha_{\pi}(i,\nu)^{k_1}\overline{\alpha_{\pi'}(j,\omega)}^{k_2}   \nonumber  \end{equation}       
since $p$ doesn't split in $EF$ we have 
$f_{\mathfrak{P}/p},f_{\mathfrak{Q}/p}\geq 2$ and since $f_{\mathfrak{P}/p}=f_{\mathfrak{P}/\mathfrak{p}}f_p\leq f_p^2$ and$f_{\mathfrak{Q}/p}=f_p' f_{\mathfrak{Q}/ \mathfrak{q}}\leq f_{p}'^{2}$ we know that $f_p,f_p'\geq 2$ so that $p$ doesn't split in $E$ or $F$.  Hence under Conjecture 2.1 we have the above sum is bounded by 
\begin{equation} \sum_{p^{k_1f_p}=p^{k_2f_p'}\leq x}(\log p)p^{k_1f_p\theta_p+k_2f_p'\theta_p'}\ll x\Big(\sum_{n\leq x}\frac{\Lambda(n)}{n^{1+\epsilon}}\Big)\ll x\exp\{-c\sqrt{\log x}\}   \nonumber  \end{equation}   
So if the collection of twisted equivalent pairs is nonempty we have the estimate
\begin{equation}  \sum_{\underset{p^{k_1}=p^{k_2}\leq x\text{   p splits in $EF$}}{k_1\leq \min\{2/(m^2\ell+1),2/(m'^2q+1)\} }}(\log p)a_{\pi\times_{BC}\pi'}(p^{k_1f_p})=\frac{x^{1+i\tau_0}}{1+i\tau_0}+O\{x\exp(-c\sqrt{\log x})\}   \nonumber  \end{equation}  
Using Hypothesis H as before and (5.2)-(5.3) we can restrict the sum to $k_1=k_2=1$ to get 
\begin{equation}  \sum_{{p\leq x}\atop {\text{  p splits completely in $EF$}}}(\log p)a_{\pi\times_{BC}\pi'}(p)=\frac{x^{1+i\tau_0}}{1+i\tau_0}+O\{ x\exp(-c\sqrt{\log x})\}  \nonumber  \end{equation}
as desired.

\bigskip

\centerline {\sc Acknowlegments}

The authors would like to thank Professor Jianya Liu and Professor
Yangbo Ye for their constant encouragement and support.  The authors would also like to thank Professor Muthu Krishnamurthy for his helpful advice and suggestions.

\end{document}